\documentclass[a4paper,10pt]{amsart}
\usepackage[all]{xy}
\usepackage{amssymb}
\usepackage[hypertex]{hyperref}
\usepackage{epic}
\usepackage{graphicx}

\renewcommand{\O}{{\mathcal{O}}}
\newcommand{\Z}{{\mathbb{Z}}}

\newcommand{\Gal}{{ \rm Gal}}

\newtheorem{pro}{Proposition}[section]
\newtheorem{lemma}[pro]{Lemma}

\newtheorem{theorem}[pro]{Theorem}

\newtheorem{definition}[pro]{Definition}

\title{Computing polynomials of the Ramanujan $\mathbf{t_n}$ class invariants.}
\author{ Elisavet Konstantinou   \and  Aristides Kontogeorgis}

\address{
Department of Information and Communication Systems Engineering, University of the Aegean, 83200 Karlovassi, Samos, Greece.
}
\email{ekonstantinou@aegean.gr}

\address{
Department of Mathematics, University of the Aegean, 83200 Karlovassi, Samos,
Greece\\ { \texttt{\upshape http://eloris.samos.aegean.gr}}
}
\email{kontogar@aegean.gr}

\date{\today{ Mathematics Subject Classification: 11R29,33E05,11R20}}

\begin{document}
\bibliographystyle{amsplain}

\maketitle

\begin{abstract}
We compute the minimal polynomials of the Ramanujan values $t_n$, where $n\equiv 11 \mod 24$,
using Shimura reciprocity law. 
These polynomials can be used for defining the Hilbert class field of the imaginary quadratic field
$\mathbb{Q}(\sqrt{-n})$, and have much smaller coefficients than the Hilbert polynomials.
\end{abstract}

\section*{Introduction}

Ramanujan on his third notebook, pages 392 and 393 in the pagination of \cite[vol. 2]{RamNotebooks}  defined the 
values
\begin{equation} \label{tndef}
t_n:=\sqrt{3} q_n^{1/18} \frac{f(q_n^{1/3}) f(q_n^3)}{f^2(q_n)}
\end{equation}
where
\[
q_n=\exp(-\pi \sqrt{n}).
\]
The function $f$ is equal to:
\begin{equation}\label{Ramt}
f(-q):=\prod_{n=1}^\infty (1-q^n) = q^{-1/24}\eta(\tau)
\end{equation}
where $q= \exp(2\pi i \tau)$, $\tau \in \mathbb{H}$ and $\eta(\tau)$ denotes the Dedekind eta-function.
Without any further explanation on how he found them,
 Ramanujan gave the following  table of polynomials $p_n(t)$ based on $t_n$ for five values of $n$:
\[
\begin{array}{|c|c|}
\hline
 n & p_n(t)\\
\hline
11 & t-1 \\
35 & t^2+t-1\\
59 & t^3+2t-1 \\
83 & t^3+2t^2+2t-1\\
107 & t^3-2t^2+4t-1\\
\hline
\end{array}
\]
In \cite{Berndt-Chan}  Bruce C.  Berndt and Heng Huat Chan 
proved that these polynomials indeed have roots
  the Ramnaujan values $t_n$. Unfortunately, their method could not be
  applied for higher values of $n$ and they asked for an efficient way of computing the polynomials $p_n$ for every $n$.
Moreover, the authors proved that if the class number of $K_n:=\mathbb{Q}(\sqrt{-n})$ is odd and 
$n \in \mathbb{N}$ is squarefree so that $n\equiv 11 \mod 24$ then
$t_n$ is a real unit generating the Hilbert class field.


It is known that the Hilbert class field can also be constructed by considering the 
irreducible polynomial of the algebraic integer $j(\theta)$ where $\theta=-1/2+i \sqrt{n}/2$.
The minimal polynomial of $j(\theta)$ is called  the Hilbert polynomial.
It is interesting to point out that the coefficients of the polynomials $p_n$ have remarkably smaller 
size compared to  the coefficients of the  Hilbert polynomials.
Therefore, finding an efficient and simple method for their construction is highly desirable and
has a direct impact to applications where the explicit construction of class fields is needed.
Problems such as primality testing/proving \cite{AtkinMorain}, the generation of elliptic curve parameters \cite{KKSZ04_icisc} and the representability of 
primes by quadratic forms \cite{Cox} can be considerably improved if the polynomials $p_n$ could be
constructed with an efficient and easily implemented way.

An explicit construction of the Hilbert class field 
has been done 
by   N. Yui  and D. Zagier in \cite{YuiZagier} using the Weber functions.
Yui and Zagier use a clever construction of a function  on quadratic forms  $ax^2 +bxy+cy^2$
that does not depend on the equivalence class of quadratic forms. The construction of 
a similar function in the case of $p_n$ polynomials seems very complicated and it is clear
that a different approach must be followed. 
Our construction came from the enforcement of Shimura reciprocity law on the values $t_n$. 
Shimura reciprocity law  has been proven to be a very powerful tool for 
attacking similar problems \cite{GeeHuatTan}, \cite{GeeBordeaux}, \cite{GeeStevenhagen} and can provide
methods for systematically determining the instances when a given function yields a class invariant
and for computing the minimum polynomial of a class invariant.


The contribution of this paper is twofold. Firstly, we prove that the values 
$t_n$ constitute class invariants for all values of $n \equiv 11 \mod 24$.
Expanding the theorem in \cite{Berndt-Chan}, we show that $t_n$ generates  the Hilbert class field
not only in the case that the class number of $K_n$ is odd but also when it is even.
Secondly, we provide an efficient method for constructing the 
irreducible polynomials $p_n$ from the  Ramanujan values $t_n$ and thus answer to the demand placed in
\cite{Berndt-Chan} for a direct and easily applicable construction method.  
Moreover, we have  implemented our method in  {\em gp-pari} \cite{pari}  and  
we present  all polynomials $p_n$ 
for all 
integers $107<n \leq  1000$, where $n\equiv 11 \mod 24$ (table 1).

The rest of the paper is organized as follows. In  the first section we fix the notation and we give the Ramanujan 
$t_n$ values in terms of the Dedekind eta function. Next, we define  six modular functions of level $72$
 $R,R_1,R_2,R_3,R_4,R_5$  and we 
compute the action of the generators of the group $SL_2(\mathbb{Z})$ on them. In section 2  we   
 prove that $t_n$ is indeed a class invariant, for all values $n\equiv 11 \mod 24$,  and in final section we employ 
Shimura reciprocity law in order to compute the conjugates of $t_n$ under the action of the 
class group and compute 
 the minimal polynomial of $t_n$.

\section{Notation}
%
%


Let $SL_2(\mathbb{Z})$ be the group of matrices with integer entries and of determinant one.
It is known \cite[cor. 1.6]{SilvII} that the group $SL_2(\mathbb{Z})$ is generated by the matrices
\[
S=\begin{pmatrix} 0 & -1 \\ 1 & 0 \end{pmatrix} \mbox{ and } 
T=\begin{pmatrix} 1 & 1 \\ 0  & 1 \end{pmatrix}.
\]
Every matrix $\begin{pmatrix} a & b \\ c & d \end{pmatrix}$ of $SL_2(\mathbb{Z})$ induces an action on the upper half plane
\[
\mathbb{H}:=\{z \in \mathbb{C}: \mathrm{Im}(z) >0 \},
\]
by sending $z \mapsto \frac{az+b}{cz+d}$.

Let $\eta$ denote the Dedekind function:
\begin{equation} \label{etadef}
\eta(\tau)=\exp(2\pi i \tau/24) \prod_{n=1}^\infty (1- q^n), \mbox{ where }  \tau\in \mathbb{H} \mbox{ and } q= \exp(2\pi i \tau). 
\end{equation}
Then $\eta$-function is transformed by $S$ and $T$ as follows \cite[prop. 8.3]{SilvII}
\begin{equation} \label{etatransf}
\eta(\tau+1)=e^{2\pi i/24} \eta(\tau) \mbox{ and } \eta(-\frac{1}{\tau})=\sqrt{-i \tau}\eta(\tau).
\end{equation}
From equation (\ref{etatransf}) we can compute the action of every element $g$ of $SL_2(\mathbb{Z})$ 
on the $\eta$-function, since $g$ can be written as a word in $S,T$.

We denote by $H_n$ the Hilbert field of $K_n:=\mathbb{Q}(\sqrt{-n})$, {\em i.e.,} the maximal Abelian unramified extension of $K_n$.
The extension $H_n/K_n$ is Galois with Galois group equal to the class group of fractional 
ideals modulo principal fractional ideals.
For imaginary quadratic fields the class group can be represented as the space of 
binary quadratic forms $ax^2+b xy +c y^2$ modulo an equivalence relation \cite[th. 5.30]{Cox}.
We will denote by $[a,b,c]$ the quadratic form $ax^2+b xy +c y^2$ and we will call two quadratic forms 
$[a_i,b_i,c_i]$ for $i=1,2$ equivalent if the  corresponding roots $\tau_i \in \mathbb{H}$ are in the 
same orbit of $SL_2(\mathbb{Z})$ acting on $\mathbb{H}$. Using the identification 
of equivalence classes of quadratic forms with the ideal class group we can define 
the structure of an abelian group on the set of equivalence classes of quadratic forms.

Let $\ell_0:=(1,1,\frac{1-d}{4})$ ($d=-n\equiv 1 \mod 4$,
 $n\in \mathbb{N}$) be the zero element in this group.
This element corresponds to the root 
\[
\tau_{\ell_0}=-\frac{1}{2}+i \frac{\sqrt{n}}{2}.
\]
Set 
\[
q_n=\exp(-\pi \sqrt{n})=-\exp(2\pi i \tau_{\ell_0}).
\]
Then 
\[
f(q_n)= f\big( -\exp(2\pi i \tau_{\ell_0})\big)=  \exp(2\pi i \tau_{\ell_0})^{-1/24} \eta(\tau_{\ell_0}),
\]
\[
f(q_n^3)= \exp(2\pi i \tau_{\ell_0})^{-3/24} \eta(3\tau_{\ell_0}),
\]
\[
f(q_n^{1/3})=(-1)^{1/18}\exp(2\pi i \tau_{\ell_0})^{-\frac{1}{3\cdot 24}}\eta(\frac{\tau_{\ell_0}}{3}+\frac{2}{3})
\]
Taking equation (\ref{tndef}) and all the above equations into consideration
 we arrive easily at the following:

\begin{lemma} \label{modtn}
 The Ramanujan value $t_n$ is given by 
\[
t_n=\sqrt{3} R_2(\tau_{\ell_0}),
\]
where 
\[
R_2(\tau)=\frac{\eta(3\tau)\eta(\frac{1}{3}\tau+\frac{2}{3})}{\eta^2(\tau)}
\]
\end{lemma}

Now let $N \in \mathbb{N}$ and  
$\Gamma(N)$ be the group 
\[
 \Gamma(N):=\left\{  \gamma \in SL_2(\mathbb{Z}),  \gamma \equiv 
\left(
\begin{array}{cc}
 1 & 0 \\
0 & 1 
\end{array}
\right) \mod N \right\}.
\]
The field  of modular functions of level $N$ consists of the  meromorphic functions  $g$ of the upper half plane $\mathbb{H}$ that are 
 invariant under the group $\Gamma(N)$, {\em i.e.} $g(\gamma \tau)=g(\tau)$ for every $\tau\in \mathbb{H}$ 
and $\gamma \in \Gamma(N)$. 
Every modular function is periodic with period $N$ and thus it admits a Fourier expansion of the form
\[
 g(q)=\sum_{\nu=-i}^\infty a_\nu q^{\nu},
\]
where $q=\exp(2\pi i \tau/N)$. We will limit ourselves to modular functions where  all coefficients of the 
Fourier expansions  are elements of the field $\mathbb{Q}(\zeta_{N})$. 
The  Galois group $\Gal(\mathbb{Q}(\zeta_{N})/\mathbb{Q})$ is isomorphic to the 
group  $\left(\frac{\mathbb{Z}}{N\mathbb{Z}}\right)^*$
by defining $\sigma_d(\zeta_{N})=\zeta_{N}^d$ for every $(d,N)=1$.

The action of the group $\Gal(\mathbb{Q}(\zeta_{N})/\mathbb{Q})$
 can be extended to 
the field of modular functions of level $N$  
with coefficients in $\mathbb{Q}(\zeta_{N})$, as follows:
\begin{equation} \label{mod72act}
g(q)^{\sigma_d} =\sum_{\nu=-i}^\infty \sigma_d(a_\nu)  q^{\nu}, \;\; \sigma_d\in Gal(\mathbb{Q}(\zeta_{N})/\mathbb{Q}).
\end{equation}

Moreover, the  action of an element $A\in GL_2(\mathbb{Z})$ on modular functions $g(q)$
can be expressed as
\begin{equation} \label{glact}
g^A=(g^B)^{\sigma_{\det(A)}}
\end{equation}
 where  $A=B\cdot \begin{pmatrix} 1 & 0 \\ 0 & \det(A) \end{pmatrix}$  and $B\in SL_2(\mathbb{Z})$.

\begin{lemma}
 The following functions are modular functions of level $72$. 
\[
 R(\tau)= \frac{\eta(3\tau)\eta(\tau/3)}{ \eta^2(\tau)}
\]
\[
R_1(\tau)= \frac{\eta(3\tau)\eta(\tau/3+1/3)}{  \eta^2(\tau)}
\]
\[
R_2(\tau)= \frac{\eta(3\tau)\eta(\tau/3+2/3)}{ \eta^2(\tau)}
\]
\[
R_3(\tau)= \frac{\eta(\tau/3)\eta(\tau/3+2/3)}{  \eta^2(\tau)}
\]
\[
R_4(\tau)= \frac{\eta(\tau/3)\eta(\tau/3+1/3)}{ \eta^2(\tau)}
\]
\[
R_5(\tau)=\frac{\eta(\tau/3+2/3) \eta(\tau/3+1/3) }{ \eta^2(\tau)}.
\]

Moreover, the element $\sigma_d: \zeta_{72} \mapsto \zeta_{72}^d$ for $(d,n)=1$ acts on them as follows:
\begin{equation} \label{111}
 \sigma_d(R)=R,
\end{equation}
\[
 \sigma_d(R_1)=
\left\{
\begin{array}{ll}
\zeta_{72}^{d-1} R_1 & \mbox{ if }  d\equiv 1 \mod 3\\
 \zeta_{72}^{d-2 } R_2 & \mbox{  if } d\equiv 2 \mod 3
\end{array}
\right.
\]
\[
 \sigma_d(R_2)=
\left\{
\begin{array}{ll}
\zeta_{72}^{2d-2} R_2 & \mbox{ if }  d\equiv 1 \mod 3\\
 \zeta_{72}^{2d-1 } R_1 & \mbox{  if } d\equiv 2 \mod 3
\end{array}
\right.
\]
\[
 \sigma_d(R_3)=
\left\{
\begin{array}{ll}
\zeta_{72}^{d-1} R_3 & \mbox{ if }  d\equiv 1 \mod 3\\
 \zeta_{72}^{d-2 } R_2 & \mbox{  if } d\equiv 2 \mod 3
\end{array}
\right.
\]
\[
 \sigma_d(R_4)=
\left\{
\begin{array}{ll}
\zeta_{72}^{2d-2} R_4 & \mbox{ if }  d\equiv 1 \mod 3\\
 \zeta_{72}^{2d-1} R_3 & \mbox{  if } d\equiv 2 \mod 3
\end{array}
\right.
\]
\[
 \sigma_d(R_5)=\zeta_{72}^{3d-3} R_5.
\]
\end{lemma}
\begin{proof}
 The fact that the above equations are indeed modular of level $72$ is a direct computation using the transformations 
of the $\eta$-functions under the generators $T,S$ of $SL_2(\mathbb{Z})$ given in (\ref{etatransf}).
The action of $\sigma_d$ given in (\ref{111}) is computed by considering the Fourier expansions of the 
$\eta$-factors of the functions $R_i$.
For instance let as compute the action of the element $\sigma_d$ on $R_2$. We begin by computing its action on 
$\eta(\tau/3+2/3)$:
\begin{equation} \label{llooppp}
 \eta(\tau/3+2/3)=\exp\left(\frac{2\pi i}{24}( \tau/3+2/3)\right) \sum_{\nu=0}^\infty a_n \exp(\frac{2 \pi i \nu }{3} \tau 
+ \frac{2 \pi i \nu}{3})= 
\end{equation}
\[
 =\exp\left(\frac{2\pi i}{24}( \tau/3)\right) \zeta_{72}^2 \sum_{\nu=0}^\infty \zeta_3^{2\nu} a_n \exp(\frac{2 \pi i \nu}{3} \tau),
\]
where $\zeta_3=\zeta_{72}^{24}$ is a primitive  third root of unity. The desired formulas of equations (\ref{111}) 
follow from the definition of the action of 
$\sigma_d$ on $\zeta_{72}$ and arguing as above.
For example (notice that from (\ref{etadef}) the Fourier expansion for $\eta(\tau)$ has rational coefficients so it 
is invariant under the action of $\sigma_d$),
\[
\sigma_d(R_2)=\sigma_d\left( \frac{\eta(3\tau) \eta(\tau/3+2/3)}{\eta^2(\tau)} \right)=\frac{\eta(3\tau)}{\eta^2(\tau)} 
\sigma_d(\eta(\tau/3+2/3)).
\] 
By (\ref{llooppp}) we have that 
\[
\sigma_d(\eta(\tau/3+2/3))=\exp\left(\frac{2\pi i}{24}( \tau/3)\right) \zeta_{72}^{2d} \sum_{\nu=0}^\infty \zeta_3^{2\nu d} a_n \exp(\frac{2 \pi i \nu}{3} \tau).
\]
If $d\equiv 1 \mod 3$ then $\zeta_3^{2\nu d}=\zeta_3^{2\nu}$,
thus 
\[
\sigma_d(R_2)=\zeta_{72}^{2d-2} R_2.
\]
If $d\equiv 2 \mod 3$ then $\zeta_3^{2\nu d}=\zeta_3^{ \nu}$,
and 
\[
\sigma_d(\eta(\tau/3+2/3))=\eta(\tau/3+1/3)) \zeta_{72}^{2d-1} 
\Rightarrow \sigma_d(R_2)=\zeta_{72}^{2d-1} R_1.
\]
\end{proof}
Later we will use some computer algebra programs in order to prove that $t_n$ is indeed a 
class invariant and  find the minimal polynomials of $t_n$.
For this reason it is convenient to have the actions  of the elements  $S,T,\sigma_d$ in matrix form.
Using (\ref{etatransf}) we give the following matrix action of the elements  $S,T,\sigma_d$ on the functions $R,R_1,\ldots,R_5$:
\begin{equation} \label{mataction}
\left(
 \begin{array}{c}
 R(\tau+1) \\
R_1(\tau+1)\\
R_2(\tau+1)\\
R_3(\tau+1)\\
R_4(\tau+1)\\
R_5(\tau+1)
 \end{array}
\right)
=A_T
\left(
 \begin{array}{c}
 R(\tau) \\
R_1(\tau)\\
R_2(\tau)\\
R_3(\tau)\\
R_4(\tau)\\
R_5(\tau)
 \end{array}
\right),
\left(
 \begin{array}{c}
 R(\frac{-1}{\tau}) \\
R_1(\frac{-1}{\tau})\\
R_2(\frac{-1}{\tau})\\
R_3(\frac{-1}{\tau})\\
R_4(\frac{-1}{\tau})\\
R_5(\frac{-1}{\tau})
 \end{array}
\right)
=
A_S
\left(
 \begin{array}{c}
 R(\tau) \\
R_1(\tau)\\
R_2(\tau)\\
R_3(\tau)\\
R_4(\tau)\\
R_5(\tau)
 \end{array}
\right),
\end{equation}
\[
\left(
 \begin{array}{c}
 \sigma_d R\\
\sigma_d R_1\\
\sigma_d R_2\\
\sigma_d R_3\\
\sigma_d R_4\\
\sigma_d R_5
 \end{array}
\right)=
A_{\sigma_d}
\left(
 \begin{array}{c}
 R \\
R_1\\
R_2\\
R_3\\
R_4\\
R_5
 \end{array}
\right),
\]
where 
\begin{equation} \label{ATAS}
A_T:={\begin{pmatrix}
 {0}&{\zeta_{72}^{3}}&{0}&{0}&{0}&{0}\cr
 {0}&{0}&{\zeta_{72}^{3}}&{0}&{0}&{0}\cr
 {\zeta_{72}^{6}}&{0}&{0}&{0}&{0}&{0}\cr
 {0}&{0}&{0}&{0}&{{1}\over{\zeta_{72}^{3}}}&{0}\cr
 {0}&{0}&{0}&{0}&{0}&{{1}\over{\zeta_{72}^{6}}}\cr
 {0}&{0}&{0}&{{1}\over{\zeta_{72}^{3}}}&{0}&{0}\cr
\end{pmatrix} },
\end{equation}
\[
A_S:=
{
\begin{pmatrix}
 {1}&{0}&{0}&{0}&{0}&{0}\cr
 {0}&{0}&{0}&1 \over{{\zeta_{72}^{3}}({ {-\zeta_{72}^{30} + \zeta_{72}^{6}} }})&{0}&{0}\cr
 {0}&{0}&{0}&{0}&{{{\zeta_{72}^{3}}\over{-\zeta_{72}^{30} + \zeta_{72}^{6}}}}&{0}\cr
 {0}&{ {-\zeta_{72}^{33} + \zeta_{72}^{9}} }&{0}&{0}&{0}&{0}\cr
 {0}&{0}&{{-\zeta_{72}^{30} + \zeta_{72}^{6}}\over{\zeta_{72}^{3}}}&{0}&{0}&{0}\cr
 {0}&{0}&{0}&{0}&{0}&{1}\cr
 \end{pmatrix} }.
\]
and 
\[
A_{\sigma_d}=
\begin{pmatrix}
1 & 0 & 0 & 0& 0 & 0 \\
0 & \zeta_{72}^{d-1} & 0 & 0 & 0& 0 \\
0 & 0 & \zeta_{72}^{2 d -d} & 0 & 0 & 0 \\
0 & 0 & 0 & \zeta_{72}^{2d-2} & 0 & 0 \\
0 & 0 & 0 & 0 & \zeta_{72}^{d-1} &0 \\
0 & 0 & 0 & 0 & 0 & \zeta_{72}^{3d-3}
\end{pmatrix}  \mbox{ if } d\equiv 1 \mod 3,
\]
\[
A_{\sigma_d}=
\begin{pmatrix}
1 & 0 & 0 & 0& 0 & 0 \\
0 & 0 &  \zeta_{72}^{d-2}  & 0 & 0& 0 \\
0 & \zeta_{72}^{2 d -1} &0 & 0 & 0 & 0 \\
0 & 0 & 0 & 0 & \zeta_{72}^{2d-1}  & 0 \\
0 & 0 & 0  & \zeta_{72}^{d-2} &0  &0 \\
0 & 0 & 0 & 0 & 0 & \zeta_{72}^{3d-3}
\end{pmatrix}  \mbox{ if } d\equiv 2 \mod 3,
\]

\section{Shimura reciprocity law and $t_n$ class invariants}
%
%


Let $K_n=\mathbb{Q}(\sqrt{-n})$ be an imaginary quadratic number field, 
$\mathcal{O}$ the ring of integers of $K_n$ and
$\theta=\frac{1}{2} - i \frac{\sqrt{n}}{2}$.
It is known that  $j(\theta)$ is an algebraic  integer that  generates the Hilbert class field of $K_n$,
and moreover that the conjugates  $j(\theta)$ under the action of the class group are  given by 
\[
  j(\theta)^{[a,-b,c]}= j(\tau_{[a,b,c]}),
\]
where $\tau_{[a,b,c]}$ is the unique root of $ax^2+bx+c$ with positive imaginary part,
{\em i.e.} $\tau_{[a,b,c]}=\frac{-b + i\sqrt{-D}}{2a}$.

There is an efficient algorithm for computing a set of non-equivalent quadratic forms and 
the minimal polynomial $f_D(x)$ of $j(\theta)$ can easily be computed from the floating point approximation of the values 
$j(\tau_{[a,b,c]})$:
\[
 f_D(x) =\prod_{[a,b,c] \in Cl(\O)} (x-j(\tau_[a,b,c])).
\]
For instance, the polynomial for the quadratic extension of discriminant $-107$ is:
 \[
f_{-107}(x)=x^3 + 129783279616\cdot 10^3 x^2 -6764523159552\cdot 10^6 x + 337618789203968\cdot 10^9.
\]
The disadvantage of using the above polynomials for the construction of the Hilbert class field is the very large size of their 
coefficients compared to the coefficients of the polynomials $p_n$. Notice that the polynomial $p_n$ for $n=107$
is equal to $x^3-2x^2+4x-1$.
In the literature there are alternative 
explicit constructions of the Hilbert class fields based on the Weber functions \cite{YuiZagier}, or 
by other modular functions \cite{GeeHuatTan}, \cite{GeeBordeaux}, \cite{GeeStevenhagen}.
In \cite{Berndt-Chan} the authors proved that also the values $t_n$ can generate the Hilbert class
field providing the following theorem.

\begin{theorem} \label{odd}
 If $n\equiv 11 \mod 24$ and the class group of $K_n$ is odd then $t_n$ generates the Hilbert field.
\end{theorem}
\begin{proof}
\cite[th. 4.1]{Berndt-Chan}.
\end{proof}

The Shimura reciprocity law can be applied  in order to compute the minimal 
polynomial of $t_n$ and  give an alternative  proof of theorem  \ref{odd} by removing the 
odd class number requirement.

In order to prove that $t_n$ is a class invariant also for the case that the class  group of $K_n$ is even,
we will use the following construction:
\begin{theorem} \label{shimura1}
Let $\mathcal{O}=\mathbb{Z}[\theta]$ be the ring of algebraic integers of the imaginary quadratic field $K$,
and assume that $x^2+Bx+C$ is the minimal polynomial of $\theta$.
Let $N>1$ be a natural number, $x_1,\ldots,x_r$ be  generators of the abelian group $\left(\mathcal{O}/N \mathcal{O} \right)^*$ and
$\alpha_i+\beta_i\theta \in \mathcal{O}$ be a representative of the class of the generator $x_i$. 
For each representative we consider 
 the matrix:
\[A_i:=\begin{pmatrix} \alpha_i-B\beta_i & -C \beta_i \\ \beta_i & \alpha_i \end{pmatrix}.\]
If $f$ is a modular function of level $N$ and if  for all matrices $A_i$ it holds that
\begin{equation} \label{act123}
f(\theta)=f^{A_i}(\theta), \mbox{ and } \mathbb{Q}(j) \subset \mathbb{Q}(f)
\end{equation}
then $f(\theta)$ is a class invariant.
\end{theorem}
\begin{proof}
\cite[Cor. 4]{GeeBordeaux}
\end{proof}

We know by lemma \ref{modtn} that the Ramanujan invariant $t_n$ can be constructed by evaluating  the modular function $\sqrt{3}R_2$ 
of level 72 at $\tau_{\ell_0}$. 
Thus,
we begin by constructing the generators of $\left(\mathcal{O}/72 \mathcal{O}\right)^*$ as Theorem ~\ref{shimura1}
dictates. Since $n\equiv 11 \mod 24$
we can take 
\[\theta=\frac{1}{2} + \frac{1}{2} \sqrt{-n}\]
 as a generator of the ring of algebraic integers of
$K=\mathbb{Q}(\sqrt{-n})$. The minimal polynomial of $\theta$ is $x^2 -x+\frac{n+1}{4}$ and thus $B$ and $C$ in
Theorem ~\ref{shimura1} are equal to $-1$ and $\frac{n+1}{4}$ respectively.
The form of the minimal polynomial implies that the prime $p=2$ stays inert in the extension $K/\mathbb{Q}$ while 
the prime $p=3$ splits.

In order to prove that $t_n:=\sqrt{3} R_2(\theta)$ is indeed a class invariant 
we have to prove that 
\begin{equation} \label{finc}
\big(\sqrt{3} R_2 \big) ^{A_i} = \sqrt{3} R_2, \mbox{ for all matrices } A_i \mbox{ and for } n\equiv 11 \mod 24.
\end{equation}
We observe that the structure  of the group $\left( \frac{\mathcal{O}}{72\mathcal{O}} \right)^*$
depends only on the value of $n \mod 72$, and there are exactly three equivalence classes 
$n \mod 72$ so that $n\equiv 11 \mod 24$, namely $n=11,35,59 \mod 72$.
Thus, we have reduced (\ref{finc}) to be checked for finite number of $n$. 

Using Chinese remainder theorem we can express the group  $\left( \frac{\mathcal{O}}{72\mathcal{O}} \right)^*$  as a direct product
\[
\left( \frac{\mathcal{O}}{72\mathcal{O}} \right)^*\cong \left(\frac{\mathcal{O}}{9\mathcal{O}}\right)^* \times 
\left(\frac{\mathcal{O}}{8\mathcal{O}}\right)^*.
\]
We will study the structure of the above two summands separately.


We  compute that 
 \[\left(\frac{\mathcal{O}}{9\mathcal{O}}\right)^*\cong \frac{\mathbb{Z}}{6\mathbb{Z}} \times \frac{\mathbb{Z}}{6 \mathbb{Z}}.\]
A selection of generators for this group is cumbersome to do by hand. We have used  brute force method, 
(we have checked all ellements one by one if there are invertible and then we have computed their orders)
   using   {\em  magma}  \cite{magma} algebra system in 
order to  compute that for $C=\frac{n+1}{4} \in \{3,9,15\}$ a set of generators is given by $7\theta+4,5$.
Moreover the group
\[
\left(\frac{\mathcal{O}}{8\mathcal{O}} \right)^* \cong \frac{\mathbb{Z}}{12 \mathbb{Z}} \times  \frac{\mathbb{Z}}{2 \mathbb{Z}}
\times  \frac{\mathbb{Z}}{2 \mathbb{Z}},
\]
and we have computed again using {\em magma} the following  selection of generators 
\[
\begin{array}{|l|l|l|}
\hline
n & \frac{n+1}{4} & \mbox{Generators} \\
\hline
11 & 3        &    \theta,7,4\theta+7       \\
35 &  9        &    5\theta+6,7,4\theta+7       \\
59 &  15      &     \theta,7,4\theta+7        \\
\hline
\end{array}
\]
From the above generators and from the Chinese remainder theorem we can construct generators for the group $\left(\mathcal{O}/72 \mathcal{O}\right)^*$ and map them to the matrices $A_i$ defined 
in theorem \ref{shimura1}. We have totally $5$ generators for the group  $\left( \frac{\mathcal{O}}{72\mathcal{O}} \right)^*$:  the first two 
are generators of the group $(\mathcal{O}/9\mathcal{O})^*$ and the last three are generators of the 
group $(\mathcal{O}/8\mathcal{O})^*$. 
In order to compute the term $f^{A_i}$ in (\ref{act123}) we have to consider any lift of $A_i$ in $GL_2(\mathbb{Z})$ and
write it as a product $w_i(S,T) \mathrm{diag}(1,\det(A_i))$, where $w_i(S,T)$ is a word in $S,T$.

The following lemma gives us the decomposition of a matrix in $SL_2(\mathbb{Z}/p^r \Z)$ as a word in the 
generators of the group  $SL_2(\mathbb{Z}/p^r \Z)$. Therefore,
we must  consider the matrices  $A_i$  modulo $8$ or $9$ and   define
$A_{i,8} \in GL_2(\mathbb{Z}/8 \mathbb{Z})$ and $A_{i,9}  \in GL_2(\mathbb{Z}/9 \mathbb{Z})$, 
so that $A_i \equiv A_{i,8} \mod 8$ and $A_i \equiv A_{i,9} \mod 9$.
 Notice that $A_{i,8}\equiv \mathrm{Id} \mod 8$ for $i=1,2$, {\em i.e.} the first two generators
and $A_{i,9} \equiv \mathrm{Id} \mod 9$ for $i=3,4,5$,  {\em i.e.} the last three generators.

\begin{lemma} \label{GeeDecompose}
Let $p^r$ be a prime power and let $\begin{pmatrix} a & b \\ c & d \end{pmatrix} \in SL_2(\mathbb{Z}/p^r \Z)$
so that either $a$ or $c$ is invertible modulo $p^r$. 
Let $\bar{S}_{p^r}=\begin{pmatrix} 0 & 1 \\ -1 & 0 \end{pmatrix}$, 
$\bar{T}_{p^r}=\begin{pmatrix} 1 & 1 \\ 0 & 1 \end{pmatrix}$ be two generators of the 
group $SL_2(\Z/p^r \Z)$.
Set $y= (1+a) c^{-1} \mod p^r$ if $(c,p)=1$, 
otherwise set $z= (1+c)a^{-1} \mod p^r$. Then 
\[
\begin{pmatrix} a & b \\ c & d \end{pmatrix} \equiv
\left\{
\begin{array}{ll}
\bar{T}_{p^r}^y \bar{S}_{p^r}  \bar{T}_{p^r} ^c \bar{S}_{p^r}  \bar{T}_{p^r}^{dy-b} \mod p^r  & \mbox{ if }  (c,p)=1 \\
\bar{S}_{p^r} \bar{T}_{p^r}^{-z}  \bar{S}_{p^r}  \bar{T}_{p^r}^{-a}  \bar{S}_{p^r}  \bar{T}_{p^r}^{bz-d} \mod p^r & \mbox{ if } (a,p)=1.
\end{array}
\right.
\]
\end{lemma}
\begin{proof}
\cite[lemma 6]{GeeBordeaux}.
\end{proof}

The generators $\bar{S}_{p^r},\bar{T}_{p^r}$ modulo $p^r=8,9$ are then lifted to elements 
$S_8,S_9,T_8,T_9 \in SL_2(\mathbb{Z})$ so that 
\[
S_8\equiv \begin{pmatrix} 0 & 1 \\ -1 & 0 \end{pmatrix} \mod 8, \mbox{ and } 
S_8 \equiv \mathrm{Id} \mod 9
\] 
\[
S_9\equiv \begin{pmatrix} 0 & 1 \\ -1 & 0 \end{pmatrix} \mod 9, \mbox{ and } 
S_9 \equiv \mathrm{Id} \mod 8
\] 
\[
T_8\equiv \begin{pmatrix} 1 & 1 \\ 0 & 1 \end{pmatrix} \mod 8, \mbox{ and } 
T_8 \equiv \mathrm{Id} \mod 9
\] 
\[
T_9\equiv \begin{pmatrix} 1 & 1 \\ 0 & 1 \end{pmatrix} \mod 9, \mbox{ and } 
T_9 \equiv \mathrm{Id} \mod 8.
\] 
Using Chinese remainder theorem we compute that 
\begin{equation} \label{liftST}
\begin{array}{lll}
S_8=T^{-1}ST^{-10}ST^{-1}ST^{-162}, & &
T_8=T^9,\\
 S_9=T^{-1}ST^{-65}ST^{-1}ST^{1096} & \mbox{ and } &
T_9=T^{-8}\\
\end{array}
\end{equation}
We observe that the elements $S_{8},T_8$ commute with $S_9,T_9$ modulo $72$.

The matrices $A_{i,p^r}$, $p^r=8$ or $9$  can be   decomposed as 
 products 
\[
A_{i,p^r} = B_{i,p^r} \begin{pmatrix} 1 & 0 \\ 0  & \det(A_{i,p^r}) \end{pmatrix},
\] 
where the matrices  $B_{i,p^r}$ have  determinant $1 \mod p^r$ and can be expressed, 
using lemma \ref{GeeDecompose}, as words
$w_{p^r}(S,T)$ in the generators $S,T$. The matrices of the form  $A_{i}$ act on 
the field of modular functions of level $72$ with coefficients in $\mathbb{Q}(\zeta_{72})$ as 
the product of $w_{8}(S,T)\cdot w_{9}(S,T)\cdot \mathrm{diag}(1,d_i)$, where $d_i$ is the determinant 
of $A_i$, {\em i.e.} the unique integer so that $d_i \equiv \det(A_{i,9}) \mod 9$ and 
$d_i \equiv \det(A_{i,8}) \mod 8$.




For instance the generator $7\theta+4$ of $\left( \mathcal{O}/9 \mathcal{O} \right)^*$  for $C=3$
corresponds to the matrix $A:=\begin{pmatrix} 11 & -21 \\ 7 & 4  \end{pmatrix}\equiv 
\begin{pmatrix} 2  & 6  \\ 7 & 4  \end{pmatrix} \mod 9$ (and to the identity matrix modulo 8). 
This is a matrix of determinant $2 \mod 9$ and it is decomposed as 
\[
\begin{pmatrix} 2  & 6  \\ 7 & 4  \end{pmatrix}  =
\begin{pmatrix} 2 &  3  \\ 7 & 2  \end{pmatrix}\begin{pmatrix} 1  &  0  \\ 0 & 2  \end{pmatrix},
\]
where $B:=\begin{pmatrix} 2 &  3  \\ 7 & 2  \end{pmatrix}$ is a matrix of determinant $1 \mod 9$.
 Using \ref{GeeDecompose} we 
find that
\[
B=\begin{pmatrix} 2  & 3  \\ 7 & 2  \end{pmatrix}=\bar{T}_9^3\bar{S}_9 \bar{T}_9^7\bar{S}_9\bar{T}_9^3.
\]
A lift of the  elements $\bar{S}_9,\bar{T}_9$  in $SL_2(\mathbb{Z})$ is given by equation (\ref{liftST}).
This means that we replace each $\bar{S}_9$ of the above formula to  $S_9=T^{-1}ST^{-65} ST^{-1}ST^{1096}$
 and each $\bar{T}_9$ to $T^{-9}$. This gives us the desired lift of $B$ to an element in $SL_2(\Z)$.
Using this lift, and the transformation matrices $A_S,A_T$ given in (\ref{ATAS})
 we compute that the  action of $A$ on the modular functions $R_i$ is given in terms of the 
following matrix:
\[
E:=
\begin{pmatrix}
0 &  0 &  0 &  \frac{-2\zeta_{72}^{18}}{3} + \frac{\zeta_{72}^6}{3} &  0 &  0 \\
0  &0 & \zeta_{72}^{15} - \zeta_{72}^3 &  0 &  0 &  0 \\
0 & 0 & 0 &  0 &  0 &  \frac{\zeta_{72}^{15}}{3} + \frac{\zeta_{72}^3}{3} \\
0 &  0 &  0 &  0 &  -\zeta_{72} ^9 &  0 \\
-2\zeta_{72}^{21} + \zeta_{72}^9 &  0 &  0 &  0  &  0 &  0 \\
0 & \zeta_{72}^{18} + \zeta_{72}^6 &  0 &  0 &  0 &  0
\end{pmatrix}.
\]
Let $V$ be the $\mathbb{Q}(\zeta_{72})$-vector space of modular functions 
generated by the elements $R,R_1,R_2,R_3,R_4,R_5$.
The vector space $V$ can be identified to the vector space $\mathbb{Q}(\zeta_{72})^6$, 
in terms of the map
\[
V \rightarrow \mathbb{Q}(\zeta_{72})^6
\]
\[
a_0 R +a_1 R_1 +\cdots +a_5 R_5 \mapsto (a_0,a_1,\ldots,a_5),
\]
$a_i \in Hom(V,\mathbb{Q}(\zeta_{72}))$.
The space $V^*=\mathbb{Q}(\zeta_{72})^6$ is the dual space of $V$ and we have to 
see how the action of the elements $T,S,\sigma_d$ act on $V^*$.
The elements $A_T,A_S$ defined in (\ref{ATAS}) act on $\mathbb{Q}(\zeta_{72})^6$ 
in terms of the transpose matrices $A_T^t,A_S^t$, while the action of $A_{\sigma_d}$ 
on $R,R_1,\ldots,R_5$ given on (\ref{ATAS}) acts on  $\mathbb{Q}(\zeta_{72})^6$ in terms of the 
contragredient action, {\em i.e.} by considering the transpose of the matrix 
$
A_{\sigma_{-d}}
$.
By chinese remainder theorem we compute that the element $\begin{pmatrix} 1  &  0  \\ 0 & 2  \end{pmatrix}$
acts on $\mathbb{Q}(\zeta_{72})$ as the automorphism $\sigma_{65}:\zeta_{72} \mapsto \zeta_{72}^{65}$.
Indeed, $65$ is an integer   $65 \equiv 2 \mod 9$ and $65 \equiv 1 \mod 8$.

Since $d=65\equiv 2 \mod 3$ we compute that the vector 
$(0,0,1,0,0,0)^t$ corresponding to the element $R_2$ is mapped to 
$A_{\sigma_{-65}} E^{\sigma_{-65}} (0,0,1,0,0,0)^t=(0,0,-1,0,0,0)$, 
where by $E^{\sigma_{-65}}$ we denote the matrix where all elements are acted on by $\sigma_{-d}$.
Notice that 
$\sqrt{3}=\zeta_{72}^{6} -\zeta_{72}^{30}.$  Indeed, the value  
 $i \sqrt{3}$ can be expressed  as a difference of two primitive $3$-roots of unity $\zeta_3,\zeta_3^2$ since 
 $i=\zeta_{72}^{18}$  and $\zeta_3=\zeta_{72}^{24}$.
Moreover, $\sigma_{-65}(\sqrt{3})=\zeta_{72}^{-6\cdot 65} -\zeta_{72}^{-30\cdot 65}=-\sqrt{3}$ and thus 
$\sqrt{3}R_2$ is left invariant.
 Following the same procedure it can be proven that $\sqrt{3} R_2$ stays invariant for all matrices $A_i$.

\begin{theorem}
The Ramanujan value $t_n$ is a class invariant for $n\equiv 11 \mod 24$.
\end{theorem}
\begin{proof}
The condition $\mathbb{Q}(j) \subset \mathbb{Q}(f)$ of theorem \ref{shimura1} is known \cite[proof of th. 4.1]{Berndt-Chan}.
By machine computation\footnote{ 
The magma program used for this computation is available in the following 
web location
\texttt{\upshape http://eloris.samos.aegean.gr/papers.html}
},
it turns out that (\ref{finc}) holds for all matrices $A_i$ and thus   $t_n$ is  a class invariant for all $n\equiv 11 \mod 24$.
\end{proof}

\section{Computing the polynomials $p_n$}

In this section we provide a method for the construction of 
 the minimal polynomial of $t_n$.
Following the article of A. Gee \cite[eq. 17]{GeeBordeaux} we give the following definition:
\begin{definition}
Let $N \in \mathbb{N}$ and 
  $[a,b,c]$ be a representative of the equivalence class of an element in the class group.
Let $p$ be a prime number and  $p^r$ be the maximum power of $p$  that divides $N$.
Assume that the discriminant $D=b^2-4ac \equiv 1 \mod 4$.
We define the matrix
\[
 A_{[a,b,c],p^r}=
\left\{
\begin{array}{ll}
\left(
\begin{array}{cc}
 a & \frac{b-1}{2} \\
0 & 1 
\end{array}
\right)
& 
\mbox{ if } p\nmid a \\
\left(
\begin{array}{cc}
 \frac{-b-1}{2} & -c \\
1 & 0 
\end{array}
\right)
& 
\mbox{ if } p\mid a \mbox{ and } p\nmid c\\
\left(
\begin{array}{cc}
 \frac{-b-1}{2}-a  & \frac{1-b}{2}-c  \\
1 & -1 
\end{array}
\right)
& 
\mbox{ if } p\mid a \mbox{ and } p \mid c.

\end{array}
\right.
\]
Chinese remainder theorem implies that 
\[
 GL_2(\mathbb{Z}/N\mathbb{Z}) \cong \prod_{p \mid N} GL_2(\mathbb{Z}/p^r \mathbb{Z}).
\]
We define $A_{[a,b,c]}$ as  the unique element in $GL_2(\mathbb{Z}/N\mathbb{Z})$ that it is mapped 
to $A_{[a,b,c],p^r}$  modulo $p^r$.
This matrix $A_{[a,b,c]}$ can be written  uniquely as a product 
\begin{equation} \label{prodexp}
A_{[a,b,c]}=
B_{[a,b,c]}
\left( 
\begin{array}{cc}
1 & 0 \\
0 & d_{[a,b,c]}
\end{array}
\right),
\end{equation}
where $d_{[a,b,c]}=\det A_{[a,b,c]}$ and $B_{[a,b,c]}$ is a matrix with determinant $1$.

\end{definition}

Shimura reciprocity law gives us \cite[lemma 20]{GeeBordeaux}  the action of 
 $[a,b,c] $    on $\sqrt{3}R_2(\theta)$ for $\theta=1/2-i \sqrt{n}/{2}$: 
\[
  \big(\sqrt{3} R_2(\theta)\big) ^{[a,-b,c]}= (\zeta_{72}^{6d_{[a,b,c]}} -\zeta_{72}^{30d_{[a,b,c]}} )
R_2\left( 
\frac{\alpha_{[a,b,c]} \tau_{[a,b,c]}+ \beta_{[a,b,c]} }
{\gamma_{[a,b,c]} \tau_{[a,b,c]} + \delta_{[a,b,c]} }
\right)^{\sigma_{d_{[a,b,c]}}},
\]
where
 $\begin{pmatrix} \alpha_{[a,b,c]} & \beta_{[a,b,c]} \\ \gamma_{[a,b,c]} & \delta_{[a,b,c]} \end{pmatrix}=A_{[a,b,c]}$
and $\tau_{[a,b,c]}$ is the (complex) root of $az^2 +bz+c$ with positive imaginary part.

If we try to implement this method in order to compute the polynomial for $t_n$ we face a 
problem. Even though we can compute a floating point approximation of the conjugate $R_2(A_{[a,b,c]})$,
it is not possible to use this approximation 
in order  to compute the action of $\sigma_{d_{[a,b,c]}}$ on it. There is however a simple approach that we 
can follow and solve this problem. 
We can express the matrix $A_{[a,b,c]}$ as a product of a matrix $B_{[a,b,c]}$ as in (\ref{prodexp}) and then
compute the expansion of  $B_{[a,b,c]}$  as a word of the matrices $S,T$.

We begin our computation by computing a full set of representatives of 
equivalence classes  $[a,b,c]$.
Since $72=2^3\cdot 3^2$ we have to compute matrices 
\[A_{[a,b,c],p^r} \in GL_2(\mathbb{Z}/p^r \mathbb{Z}) \mbox{ for } p^r=8,9.\]
Then we compute the determinant $d_{[a,b,c],p^r}$  
of the matrix $A_{[a,b,c],p^r}$ 
and we find a decomposition 
\[
A_{[a,b,c],p^r} =B_{[a,b,c],p^r} \begin{pmatrix}
           1 & 0 \\
           0 & d_{[a,b,c],p^r}
                         \end{pmatrix}.
\]
The matrices  $B_{[a,b,c],p^r}$ are elements of $SL_2(\mathbb{Z}/ p^r \mathbb{Z})$ and can 
be    written as words of $\bar{S}_{p^r},\bar{T}_{p^r}$ 
by using lemma \ref{GeeDecompose}.

So, if  $B_{[a,b,c],8}=w(\bar{T}_8,\bar{S}_8)$ and $B_{[a,b,c],9}=w'(\bar{T}_9,\bar{S}_9)$ 
are the decompositions of $B_{[a,b,c],p^r}$ as words of $\bar{S}_{p^r},\bar{T}_{p^r}$ we 
take the lift
\[
B_{[a,b,c]}=w(T_8,S_8) w'(T_9,S_9) \in SL_2(\Z),
\]
and the corresponding action on the functions $R,R_1,\ldots,R_5$ is computed by using  (\ref{mataction}).

The determinants $d_{[a,b,c],8} \in \Z/8\Z$ and $d_{[a,b,c],9} \in \Z/9\Z$ can be lifted in an element $d_{[a,b,c]}$
of $\Z/ 72 \Z$ (so that it reduces to $d_{[a,b,c],8} \mod 8$ and  $d_{[a,b,c],9} \mod 9$ respectively) by using 
again the Chinese remainder theorem.

The desired polynomial $p_n$
can then be computed:
\[
p_n(t) =\prod_{[a,b,c]} \left(t-\big(\sqrt{3} R_2(\frac{-1+i \sqrt{n}}{2} ) \big)^{[a,-b,c]} \right).
\]
We have used the {\em gp-pari}  \footnote{The pari program used for this computation is 
available at { \texttt{\upshape http://eloris.samos.aegean.gr/papers.html}}. } program in order to perform this computation.
The resulting polynomials $p_n$ for $107 \leq n < 1000$ are given in table 1.

\begin{table}[H] \label{pnpoly1}
\caption{Polynomials $p_n$ for $107 \leq n < 1000$.}

\[
\begin{array}{|l|l|}
\hline
n & p_n(t)\\
\hline
107 &
{x^{3} - {{2}}x^{2} + {{4}}x - 1} \\
131 &
{x^{5} + x^{4} - x^{3} - {{3}}x^{2} + {{5}}x - 1} \\
155 &
{x^{4} + {{2}}x^{3} + {{5}}x^{2} + {{4}}x - 1} \\
179 &
{x^{5} - {{2}}x^{4} + {{5}}x^{3} - x^{2} + {{6}}x - 1} \\
203 &
{x^{4} - {{3}}x^{3} + {{7}}x - 1} \\
227 &
{x^{5} - {{5}}x^{4} + {{9}}x^{3} - {{9}}x^{2} + {{9}}x - 1} \\
251 &
{x^{7} + {{5}}x^{6} + {{6}}x^{5} - {{2}}x^{4} - {{4}}x^{3} + {{2}}x^{2} + {{9}}x - 1} \\
275 &
{x^{4} - x^{3} + {{6}}x^{2} - {{11}}x + 1} \\
299 &
{x^{8} + x^{7} - x^{6} - {{12}}x^{5} + {{16}}x^{4} - {{12}}x^{3} + {{15}}x^{2} - {{13}}x + 1}\\
323 &
{x^{4} - x^{3} + {{4}}x^{2} + {{13}}x - 1}\\
347 &
{x^{5} + {{7}}x^{4} + {{21}}x^{3} + {{27}}x^{2} + {{13}}x - 1} \\
371 &
{x^{8} + {{9}}x^{6} - {{10}}x^{5} + {{14}}x^{4} + {{8}}x^{3} - {{23}}x^{2} + {{18}}x - 1} \\
395 &
{x^{8} - x^{7} + {{5}}x^{6} + {{16}}x^{5} + {{28}}x^{4} + {{24}}x^{3} + {{27}}x^{2} + {{17}}x - 1} \\
419 &
{x^{9} - {{6}}x^{8} + {{12}}x^{7} - {{7}}x^{6} + {{12}}x^{5} - {{8}}x^{4} + {{31}}x^{3} + {{10}}x^{2} + {{20}}x - 1}\\
443 &
{x^{5} - {{4}}x^{4} - {{3}}x^{3} + {{17}}x^{2} + {{22}}x - 1} \\
467 &
{x^{7} + {{6}}x^{6} + {{7}}x^{5} - {{3}}x^{4} + {{3}}x^{3} - {{23}}x^{2} + {{26}}x - 1}\\
491 &
{x^{9} + x^{8} + {{16}}x^{7} + {{2}}x^{6} + {{37}}x^{5} - {{31}}x^{4} + {{44}}x^{3} - {{40}}x^{2} + {{29}}x - 1}\\
515 &
{x^{6} + {{8}}x^{5} + {{32}}x^{4} + {{60}}x^{3} + {{68}}x^{2} + {{28}}x - 1}\\
539 &
{x^{8} - {{6}}x^{7} + {{28}}x^{6} - {{56}}x^{5} + {{77}}x^{4} - {{56}}x^{3} + {{28}}x^{2} - {{34}}x + 1} \\
563 &
{x^{9} + {{4}}x^{8} + {{6}}x^{7} - {{11}}x^{6} + {{44}}x^{5} - {{76}}x^{4} + {{91}}x^{3} - {{64}}x^{2} + {{38}}x - 1} \\
587 &
{x^{7} + x^{6} + {{16}}x^{5} - {{12}}x^{4} + {{20}}x^{3} + {{24}}x^{2} + {{39}}x - 1} \\
611 &
{x^{10} - {{8}}x^{9} + {{35}}x^{8} - {{62}}x^{7} - x^{6} + {{116}}x^{5} - {{65}}x^{4} - {{100}}x^{3} + {{125}}x^{2} - {{46}}x + 1} \\
635 &
{x^{10} - {{11}}x^{9} + {{50}}x^{8} - {{121}}x^{7} + {{201}}x^{6} - {{192}}x^{5} + {{87}}x^{4} + {{51}}x^{3} - {{98}}x^{2} + {{49}}x - 1} \\
659 &
{x^{11} - {{7}}x^{10} + {{7}}x^{9} + {{27}}x^{8} + {{19}}x^{7} - {{43}}x^{6} - {{5}}x^{5} + {{91}}x^{4} + {{157}}x^{3} + {{97}}x^{2} + {{49}}x - 1} \\
683 &
{x^{5} + {{6}}x^{4} - {{5}}x^{3} - {{41}}x^{2} + {{56}}x - 1} \\
707 & 
{x^{6} + {{4}}x^{5} + {{30}}x^{4} + {{72}}x^{3} + {{108}}x^{2} + {{58}}x - 1} \\
731 & 
           x^{12} + {{7}}x^{11} + {{25}}x^{10} + {{12}}x^{9} + {{41}}x^{8} + {{9}}x^{7} + \\
        & +{{92}}x^{6} + {{73}}x^{5} - {{133}}x^{4} + {{216}}x^{3} - {{153}}x^{2} + {{67}}x - 1 \\
755 &
         x^{12} - {{2}}x^{11} + {{18}}x^{10} + {{50}}x^{9} + {{82}}x^{8} + {{182}}x^{7} + {{360}}x^{6} + {{522}}x^{5} + \\
       & +  {{598}}x^{4} + {{486}}x^{3} + {{262}}x^{2} + {{66}}x - 1 \\
779 &
{x^{10} + {{8}}x^{9} + {{24}}x^{8} - {{8}}x^{7} - {{11}}x^{6} + {{26}}x^{5} + {{81}}x^{4} + {{220}}x^{3} + {{98}}x^{2} + {{74}}x - 1}\\
803 &
{x^{10} + {{3}}x^{9} + {{26}}x^{8} + {{11}}x^{7} - {{65}}x^{6} + {{16}}x^{5} + {{7}}x^{4} - {{83}}x^{3} + {{150}}x^{2} - {{83}}x + 1} \\
827 &
{x^{7} - {{7}}x^{6} + {{38}}x^{5} - {{54}}x^{4} + {{112}}x^{3} - {{146}}x^{2} + {{89}}x - 1} \\
851 &
{x^{10} - {{7}}x^{9} - x^{8} + {{86}}x^{7} + {{69}}x^{6} - {{201}}x^{5} - {{219}}x^{4} + {{94}}x^{3} + {{103}}x^{2} - {{95}}x + 1} \\
875 &
{x^{10} - {{10}}x^{9} + {{25}}x^{8} + {{10}}x^{7} + {{15}}x^{6} + {{94}}x^{5} - {{35}}x^{4} - {{120}}x^{3} + {{85}}x^{2} + {{100}}x - 1} \\
899 &
{x^{14} + {{16}}x^{13} + {{97}}x^{12} + {{308}}x^{11} + {{666}}x^{10} + {{1086}}x^{9} + } \\
       & + {{{1490}}x^{8} + {{1766}}x^{7} + {{1800}}x^{6} + {{1556}}x^{5} + {{998}}x^{4} + {{698}}x^{3} + {{229}}x^{2} + {{106}}x - 1} \\
923 & 
{x^{10} - x^{9} + {{30}}x^{8} - {{81}}x^{7} - {{29}}x^{6} + {{56}}x^{5} + {{211}}x^{4} - {{27}}x^{3} - {{110}}x^{2} - {{115}}x + 1} \\
947 &
{x^{5} + {{5}}x^{4} + {{7}}x^{3} - {{103}}x^{2} + {{125}}x - 1} \\
971 &
{x^{15} - x^{14} + {{21}}x^{13} + {{133}}x^{12} + {{264}}x^{11} + {{310}}x^{10} + {{216}}x^{9} +} \\ 
        & +{{{62}}x^{8} - {{100}}x^{7} - {{300}}x^{6} + {{152}}x^{5} + {{338}}x^{4} + {{79}}x^{3} - {{285}}x^{2} + {{135}}x - 1} \\
995 &
{x^{8} + {{12}}x^{7} + {{59}}x^{6} + {{78}}x^{5} + {{12}}x^{4} + {{66}}x^{3} + {{289}}x^{2} + {{140}}x - 1}\\
\hline
\end{array}
\]
\end{table}

\providecommand{\bysame}{\leavevmode\hbox to3em{\hrulefill}\thinspace}
\providecommand{\MR}{\relax\ifhmode\unskip\space\fi MR }
\providecommand{\MRhref}[2]{%
  \href{http://www.ams.org/mathscinet-getitem?mr=#1}{#2}
}
\providecommand{\href}[2]{#2}

\end{document}